\documentclass[12pt]{amsart}
\usepackage{amssymb}
\usepackage[curve]{xypic}
\newcommand{\comment}[1]{\marginpar{\sffamily{\raggedright\noindent\tiny #1
   \par}\normalfont}}
\newcommand{\Omit}[1]{\begin{tiny}#1\end{tiny}}
\renewcommand{\Omit}[1]{}
\renewcommand{\comment}[1]{}

\newbox\mybox
\def\overtag#1#2#3{\setbox\mybox\hbox{$#1$}\hbox to
  0pt{\vbox to 0pt{\vglue-#3\vglue-\ht\mybox\hbox to \wd\mybox
      {\hss$\ss#2$\hss}\vss}\hss}\box\mybox}
\def\undertag#1#2#3{\setbox\mybox\hbox{$#1$}\hbox to 0pt{\vbox to
    0pt{\vglue#3\vglue\ht\mybox\hbox to \wd\mybox
      {\hss$\ss#2$\hss}\vss}\hss}\box\mybox}
\def\lefttag#1#2#3{\hbox to 0pt{\vbox to 0pt{\vss\hbox to
      0pt{\hss$\ss#2$\hskip#3}\vss}}#1}
\def\righttag#1#2#3{\hbox to 0pt{\vbox to 0pt{\vss\hbox to
      0pt{\hskip#3$\ss#2$\hss}\vss}}#1}
\let\ss\scriptstyle

\def\Dot{\lower.2pc\hbox to 2.5pt{\hss$\bullet$\hss}}
\def\Circ{\lower.2pc\hbox to 2.5pt{\hss$\circ$\hss}}
\def\Vdots{\raise5pt\hbox{$\vdots$}}
\def\splicediag#1#2{\xymatrix@R=#1pt@C=#2pt@M=0pt@W=0pt@H=0pt}

\renewcommand\frame[2][3pt]{\hbox{$\vcenter{\hbox{\vrule\vbox
{\hrule\kern#1\hbox{\kern#1$#2$\kern#1}\kern#1\hrule}\vrule}}$}}

\def\disc{\det}

\newcommand\lineto{\ar@{-}}
\newcommand\dashto{\ar@{--}}
\newcommand\dotto{\ar@{.}}
\newcommand{\C}{\mathbb C}
\newcommand{\E}{\mathbb E}
\newcommand{\N}{\mathbb N}
\newcommand{\Z}{\mathbb Z}
\newcommand{\Q}{\mathbb Q}

\newtheorem*{ECtheorem}{End Curve Theorem}

\newtheorem*{theorem*}{Theorem}
\newtheorem{theorem}{Theorem}[section]
\newtheorem{proposition}[theorem]{Proposition}

\theoremstyle{definition}

\newtheorem*{examples*}{Examples}
\newtheorem{definition}[theorem]{Definition}

\newtheorem*{remark*}{Remark}
\newtheorem{conjecture}[theorem]{Conjecture}
\newtheorem*{questions*}{Questions}
\newtheorem*{Cconjecture*}{Casson Invariant Conjecture}
\newtheorem*{Mconjecture*}{Milnor Fiber Conjecture}
\begin{document}
\title[Splice diagrams and splice-quotient surface singularities]
{Splice diagrams and splice-quotient surface singularities}
 \author{Jonathan Wahl}
\address{Department of Mathematics\\The University of North
  Carolina\\Chapel Hill, NC 27599-3250} \email{jmwahl@email.unc.edu}
\maketitle

It has been my great pleasure to collaborate on many papers with Walter Neumann, a terrific mathematician and friend.   I'll try to present an overall account of some of our work, primarily on splice diagrams and splice-quotient singularities.  Our goal was to be able to write down explicit equations of certain kinds of normal surface singularities, given only the topology of their links.  I mention our motivations, and why and how we were led to certain kinds of questions, leaving detailed proofs (and careful definitions) to the original papers.  One can also consult interesting more recent work of others, such as \cite{cueto}, which uses a tropical geometry approach. 

\section{Provocative early work of Walter}
While I first became friendly with Walter in 1972 at the Institute for Advanced Study,  our collaboration did not begin until 1987, well after he had written two papers about normal surface singularities which I found very provocative. 

 Recall that if $(X,0)$ is a complex normal surface singularity, its \emph{link} $\Sigma$ is the oriented $3$-manifold which is its neighborhood boundary, i.e., $X$ intersected with a small sphere about $0$.  From the minimal good resolution $\pi:(\tilde{X},E)\rightarrow (X,0)$, where $\pi^{-1}(0)=E$ is a union of smooth curves with strong normal crossings, $\Sigma$ is the boundary of a tubular neighborhood of $E$.  One associates to E a plumbing graph $\Gamma$, encoding the genera, self-intersections, and intersections of the curves.  From $\Gamma$, one can reconstruct $\Sigma$ by the \emph{plumbing construction}, as made clear in seminal papers by Mumford and Hirzebruch.   Thus the plumbing graph $\Gamma$ determines the topology of $\Sigma$.   Walter's obviously important theorem proves that not only does $\Gamma$ determine the $3$-manifold $\Sigma$, but the converse is true as well.
\begin{theorem}\cite{neumann81} Suppose $\Sigma$ is a singularity link.  Then $\Sigma$ determines the minimal good resolution graph $\Gamma$.
\end{theorem}
     In fact, with the exception of lens spaces (from cyclic quotient singularities) and some two-torus bundles over the circle (from ``cusps''), already $\pi_1(\Sigma)$ determines $\Gamma$.
     
    The theorem was proved by a calculus for general plumbing diagrams (not necessarily negative-definite), and implied other basic results (e.g.,  $S^1\times S^2$ is not a singularity link).  Thus, any invariant concocted from $\Gamma$ is topological, and not just a singularity invariant.  On the other hand, the geometric genus and Milnor number depend in general on analytic information.
 
A later result of Walter's involved certain singularities with a rational homology sphere (or $\Q$HS) link; that is, $E$ is a tree of smooth rational curves.  Then $H_1(\Sigma; \Z)\equiv H_1$ is the finite \emph{discriminant group} $D(\Gamma)$, calculated for example as the cokernel of the intersection matrix $(E_i\cdot E_j)$; thus, $\Sigma$ has a finite \emph{universal abelian covering} (or UAC) $\Sigma'\rightarrow \Sigma$.  That covering is realized by a finite map of normal germs $(X',0)\rightarrow (X,0)$, which is a quotient by an action of the group $H_1$; this is the UAC of $(X,0)$.  

A weighted-homogeneous singularity with $\Q$HS link has resolution graph

$$
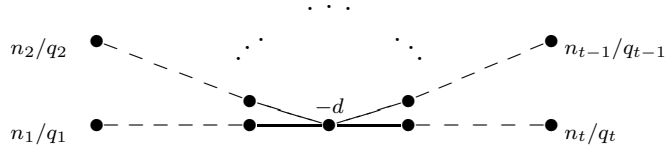
\xymatrix@R=4pt@C=24pt@M=0pt@W=0pt@H=0pt{\\
\lefttag{\bullet}{n_2/q_2}{8pt}\dashto[ddrr]&
&\hbox to 0pt{\hss\lower 4pt\hbox{.}.\,\raise3pt\hbox{.}\hss}
&\hbox to 0pt{\hss\raise15pt
\hbox{.}\,\,\raise15.7pt\hbox{.}\,\,\raise15pt\hbox{.}\hss}
&\hbox to 0pt{\hss\raise 3pt\hbox{.}\,.\lower4pt\hbox{.}\hss}
&&\righttag{\bullet}{n_{t-1}/q_{t-1}}{8pt}\dashto[ddll]\\
\\
&&\bullet\lineto[dr]&&\bullet\lineto[dl]\\
\lefttag{\bullet}{n_1/q_1}{8pt}\dashto[rr]&&
\bullet\lineto[r]&\overtag{\bullet}{-d}{8pt}\undertag{}{}{6pt}\lineto[r]&\bullet
\dashto[rr]&&\righttag{\bullet}{n_{t}/q_{t}}{8pt}\\&~\\&~\\&~\\&~}
$$
Here,  the continued fraction expansion $n/q=b_1-1/b_2-\cdots -1/b_s$ represents a string of rational curves emanating from the center:
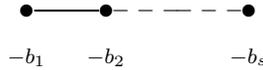
\begin{equation}
\label{step1}
\xymatrix@R=6pt@C=24pt@M=0pt@W=0pt@H=0pt{
\\&\undertag{\bullet}{-b_1}{6pt}\lineto[r]&\undertag{\bullet}{-b_2}{6pt}
\dashto[r]&\dashto[r]
&\undertag{\bullet}{-b_s}{6pt}\\
&&&&\\
&&&&\\
&&&&\\
}
\end{equation}
(We shall assume that $t\geq 3$.) The analytic type of $(X,0)$ is uniquely determined by the graph plus the location of the $t$ intersection points on the central curve; in fact, one can write $X=\text{Spec}\oplus_{i=0}^{\infty} A_i$, where each  $A_i$ is computable as a cohomology group from data on the central $\mathbb P^1$ (\cite{pinkham}, (5.1)).  Still, equations defining the singularity might be quite complicated, whereas the UAC is simple.

\begin{theorem}\cite{n83}  Notation as above, the UAC of the weighted-homogeneous singularity $(X,0)$ is a Brieskorn complete intersection $(X',0)\subset (\C^t,0)$, defined by equations $$\sum_{j=1}^t a_{ij}z_j^{n_j}=0,\ \ \ i=1,\cdots, t-2.$$  In addition,  the discriminant group $H_1$ acts on $X'$ diagonally on the $t$ coordinates $z_j$.
\end{theorem}

The matrix $(a_{ij})$ in the Theorem satisfies the \emph{Hamm condition} \cite{h}, that all maximal minors have full rank.

I found this a very surprising result!   It is easy to show that $(X',0)$ is Gorenstein; but why is it a complete intersection?  And, while its defining equations must be weighted-homogeneous, why are they of the very special Brieskorn type?  Walter's proof was simply to write down these equations, find an appropriate representation of the discriminant group $H_1$ on $\C^t$, and then to prove that (for judiciously chosen $a_{ij}$) the quotient is the singularity we started with.

Attempts to find a natural general statement which includes this Theorem led years later to our definition of splice-quotient singularities in \cite{nw1} (see Section 5).

\section{casson invariant and splice diagrams}
 In 1987 I attended a lecture given by Sir Michael Atiyah at the Weyl Centennial Conference at Duke University, entitled ``New invariants of three and four-manifolds"\cite{at}.  Andrew Casson had discovered a new integer invariant $\lambda(\Sigma)$ for a compact integral homology $3$-sphere (or $\Z$HS); its mod $2$ reduction is the Rokhlin invariant (see e.g. \cite{sav}, (12.1)).  Naturally I asked about its value for the link $\Sigma(p,q,r)$ of the Brieskorn singularity $x^p+y^q+z^r=0$, which is a $\Z$HS when $p, q, r$ are pairwise relatively prime.  Atiyah answered that based on developing work by Fintushel-Stern, $\lambda$ was a complicated expression involving the number of integer points in some three-dimensional region. I guessed it was related to the singularity itself; a call to Walter completely clarified the situation.
 
 First, $\lambda(\Sigma)$ is minus one-half the ``number'' of equivalence classes of $SU(2)$ representations of $\pi_1(\Sigma)$, where \emph{number} is an algebraic count done using a Heegaard decomposition of $\Sigma$.  Walter explained that L. Greenberg had found an actual count of such representations, given by the number of integer points in some region, and Fintushel-Stern were proving that $\lambda(\Sigma(p,q,r))$ involves simply the actual count.  But then we realized that a formula of Brieskorn \cite{br} implied that the corresponding expression for $\lambda$ was equal to $1/8$ the signature of the intersection pairing on the second homology of the Milnor fiber of the singularity.  (As the pairing is even, non-degenerate, and unimodular, a basic result implies that the signature is divisible by $8$.)  Putting together work of others, we had realized a beautiful and provocative expression for the Casson invariant of these particular links.  Our paper \cite{nw90} contained the appealing yet wildly optimistic
 
 \begin{Cconjecture*} (CIC) Suppose $(X,0)$ is the link of a hypersurface (or complete intersection) singularity whose link $\Sigma$ is an integral homology sphere.  Then the Casson invariant of $\Sigma$ equals one-eighth the signature of the Milnor fiber.
 \end{Cconjecture*}
 
 \begin{theorem} \cite{nw90}  The Casson Invariant Conjecture is true in the following cases:
\begin{enumerate}
\item The Brieskorn complete intersection $\sum_{j=1}^t a_{ij}z_j^{n_j}=0,\ \ \ i=1,\cdots, t-2$, with pairwise relatively prime $n_i$.
\item The singular germ $z^n+f(x,y)=0$, where $f(x,y)=0$ is an irreducible plane curve singularity and $n$ is relatively prime to all the numbers among the Puiseux pairs.
\item The complete intersection singularity in $\C^4$ defined by $$x^n=u^{n+1}+v^ny,\ y^n=v^{n+1}+u^nx.$$ 
\end{enumerate} 
\end{theorem}

The method was correct but mathematically disappointing: we compute and compare the Casson invariant and signature in each case.  Walter realized that the key was the use of \emph{splice diagrams}, introduced earlier by Siebenmann, and extensively studied by Eisenbud-Neumann \cite{e-n}.   \emph{Splicing} refers to a different topological construction than plumbing, allowing one to make a new $\Z$HS from two old ones.  Specifically, if $(\Sigma_i, K_i)$ is a pair of a homology $3$-sphere and a knot, one removes the interiors of closed tubular neighborhoods $K_1\times D^2$ of $K_1$ and $D^2\times K_2$ of $K_2$, and pastes the $\Sigma_i$ together along the torus boundaries $K_1\times S^1$ and $S^1\times K_2$, in other words switching the roles of meridian and longitude (and paying attention to orientations!)  This process yields a new homology sphere, which in our situation can be represented by a \emph{splice diagram}.								
															
A splice diagram $\Delta$ is a finite tree with vertices
only of valency 1 (``\emph{leaves}'') or $\ge3$ (``\emph{nodes}'') and
with a collection of integer weights at each node, associated to the
edges departing the node. The Brieskorn link $\Sigma(p,q,r)$ is represented by the diagram 
  $$
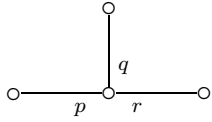
\xymatrix@R=8pt@C=30pt@M=0pt@W=0pt@H=0pt{
&&\circ\\
&&\\
&&\\
&\circ\lineto[r]_(.7)p&\circ\lineto[r]_(.3)r\lineto[uuu]_(.3)q&\circ\\\\
&\\
}
$$
The next example is the diagram from splicing $\Sigma(2,3,7)$ with $\Sigma(2,5,11)$ along the knots given by the last coordinates:
$$\splicediag{12}{30}{
\circ&&&\circ\\
&\circ\lineto[ul]_(.25){2}\lineto[dl]^(.25)3
&\circ\lineto[dr]_(.25){5}\lineto[ur]^(.25){2}
\lineto[l]_(.2){11}_(.8){7}\\
\circ&&&\circ
}$$
For an edge connecting two nodes in a splice diagram the \emph{edge
  determinant} is the product of the two weights on the only edge minus the
product of the weights adjacent to the edge. In the above
example, the one edge connecting two nodes has edge determinant
$77-60=17$.


\begin{theorem}\cite{e-n}
The homology spheres that are singularity links are in one-one
correspondence with splice diagrams satisfying 
\begin{itemize}
\item the weights around a node are positive and pairwise coprime;
\item the weight on an edge ending in a leaf is $>1$;
\item all edge determinants are positive.
\end{itemize}
\end{theorem}

From the resolution diagram for a homology sphere singularity link, one associates a splice diagram as follows: disregard all vertices of valency $2$; at each node and outgoing edge, take the absolute value of the determinant of the outer diagram.  For example, the splice diagram above arises from the resolution diagram
$$
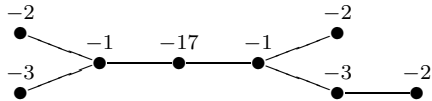
\xymatrix@R=6pt@C=24pt@M=0pt@W=0pt@H=0pt{
\\
&\overtag{\bullet}{-2}{8pt}&&&&\overtag{\bullet}{-2}{8pt}\\
{}&&\overtag{\bullet}{-1}{8pt}\lineto[ul]\lineto[dl]\lineto[r]&
\overtag{\bullet}{-17}{8pt}&\overtag{\bullet}{-1}{8pt}\lineto[ur]\lineto[dr]\lineto[l]&\\
&\overtag{\bullet}{-3}{8pt}&&&&\overtag{\bullet}{-3}{8pt}\lineto[r]&\overtag{\bullet}{-2}{8pt}}
$$

There are methods to compute the resolution diagram
from the splice diagram; see a general discussion of plumbing versus splicing in the Appendix to \cite{nw1}. 

One reason to consider splice diagrams for the CIC is that the Casson invariant is known to be additive under splicing.  All $\Z$HS singularity links admit splice diagrams built on the $\Sigma(p,q,r)$ examples, so producing these diagrams for the examples in Theorem 2.1 allows a direct calculation of the Casson invariants.  The signatures in these cases, however, are generally extremely difficult to compute, and we had to use several different methods.  

One consequence of the CIC would be that for a complete intersection singularity with $\Z$HS link, the topology of the link determines the geometric genus.  (We foolishly guessed originally that we could replace ``complete intersection'' by ``Gorenstein'', but \cite{nem} found a counterexample.)  Such a surprising result should give one pause before being too confident of its correctness.  On the other hand, the Casson invariant of a $\Z$HS  $\Sigma$ is one-eighth the signature of \emph{some} simply-connected spin-manifold whose boundary is $\Sigma$; the CIC asserts that the Milnor fiber is such a manifold.

\section{Singularities of splice type}
After the completion in 1989 of our paper on the Casson invariant, Walter (and I) pursued other topics in the 1990's, he working primarily on hyperbolic $3$-manifolds and geometric group theory.  But in 1998 Walter spent a semester at Duke University, at which time we resumed our discussion of two big questions (beyond the CIC) left open for us from \cite{nw90}. 

The first was the dearth of explicit examples of hypersurface and complete intersection singularities with $\Z$HS links (so, the intersection matrix $(E_i\cdot E_j)$ has determinant $\pm 1$).  Could one write down equations of such singularities for at least some splice diagrams from Theorem 2.2? The complete intersection examples of Theorem 2.1(3) gave a hint; we were lucky that Henry Laufer had provided equations for us from the corresponding resolution graph in the simplest case $n=2$.

Second, since the Casson invariant adds under splicing, perhaps the CIC is true because one can ``splice'' singularities and their Milnor fibers, matching  splicing on the link level, and so that on the Milnor fiber level the signatures add.  Both these issues were eventually discussed in \cite{nw1} and \cite{nw2}.

 A splice diagram with a single node and $t$ leaves arises from a Brieskorn complete intersection in $\C^{t}$ as in Theorem 2.1.1.  A variable is assigned to each of the leaves, and one defines a singularity by $t-2$ generic linear combinations of monomials from each leaf.  Weighted-homogeneity is given by assigning to each leaf a weight which is the product of all the other weights from the other edges.  In addition, note that adding to the equations terms of weight greater than or equal to the weight $n_1n_2\cdots n_t$ gives other singularities with the same resolution and splice diagrams.  
 
 Consider now the two-node splice diagram below, where we have already added a variable to each leaf:
$$\splicediag{12}{30}{
\lefttag{\circ}{X}{6pt}&&&\righttag{\circ}{Z}{6pt}\\
&\circ\lineto[ul]_(.25){p}\lineto[dl]^(.25)q
&\circ\lineto[dr]_(.25){q'}\lineto[ur]^(.25){p'}
\lineto[l]_(.2){r'}_(.8){r}\\
\lefttag{\circ}{Y}{6pt}&&&\righttag{\circ}{W}{6pt}
}$$
Recall that $p,q,p',q'$ are $\geq 2$, the triples $p,q,r$ and $p',q',r'$ are each pairwise relatively prime, and $rr'>pqp'q'$.
\begin{definition} The diagram satisfies the \emph{semigroup conditions} if $r'\in \N(p,q)$ (the semigroup generated by $p$ and $q$) and $r\in \N(p',q')$.  
\end{definition}
Since $r'\geq (p-1)(q-1)$ implies that $r'\in \N(p,q)$, the edge determinant condition does guarantee that at least one of $r,r'$ is in the appropriate semigroup.  

\begin{theorem}  For the splice diagram above, assume the semigroup conditions are satisfied, and write $$r'=\alpha p+\beta q, r=\gamma p'+\delta q'.$$   Then the isolated complete intersection singularity given by
$$X^p+Y^q+Z^{\delta}W^{\gamma}=0$$
$$Z^{p'}+W^{q'}+X^{\beta}Y^{\alpha}=0$$
has $\Z$HS link with the given splice diagram.
\end{theorem}

In analogy with the one-node case, we have first assigned a variable to each leaf.  Next, to each node we associate a weight which is the product of the weights of the surrounding edges (so, $pqr$ for the left node).  Third, fixing a node we assign a weight to each leaf, by taking the product of the weights adjacent to (but not on) the route from the node to that leaf; for the left node, this weight for the upper left leaf is $qr$, while the weight for the upper right leaf is $pqq'$.  Crucially, the semigroup condition involving $r$ allows one to write the first equation, the sum of three monomials for the three edges of the node, each of which is homogeneous with respect to these weights, of total degree $pqr$.  

Thus, with the weights from the left node, the first equation is weighted homogeneous, while the second is not.  Still,  the associated graded is an integral domain, and its normalization is the singularity $X^p+Y^q+T^r=0$.  One starts to resolve the singularity with a weighted blow-up of $\C^4$.  A key point is that setting the variable $Z$ or $W$ equal to $0$ gives a monomial curve, corresponding to the ``$r$'' knot in $\Sigma(p,q,r)$.  But adding higher weight terms to the equations would not affect the topology.

One can similarly write down a system of equations for an arbitrary splice diagram $\Delta$, as long as appropriate semigroup conditions are assumed.  Start by assigning a variable $z_w$ to every leaf $w$.  Next, for each node $v$, take as weight $d_v$ the product of the weights $d_{ve}$ on the $t$ surrounding edges.  To every leaf $w$, define its $v$-weight to be $\l_{vw}$, the product of all weights adjacent to (but not on) the path from $v$ to $w$; and define $\l'_{vw}$ to be the same product excluding the weights around $v$.  The \emph{semigroup condition} means:  the weight on any edge $e$ is in the semigroup generated by the $\l'_{vw}$ of all the outer leaves on the $e$-side of $v$.  Writing $d_{ve}=\Sigma_w\alpha_{vw}\l'_{vw}$, it follows that the \emph{admissible monomial}  $M_{ve}=\prod_w
z_w^{\alpha_{vw}}$ has weight $d_v$, the weight of the node.   Now take $t-2$ linear combinations of those $t$ monomials so that the coefficient matrix $(a_{vie})$ has all maximal minors of full rank, and add to each terms $H_{vi}$ of weight $>d_v$:
 $$\Sigma_ea_{vie}M_{ve}+H_{vi}=0, \ i=1,\cdots t-2.$$
 Finally, choose such equations for every node $v$ of $\Delta$.  The following is a special case of Theorem 5.4 below, proved originally in \cite{nw1}.

\begin{theorem}  For any splice diagram $\Delta$ satisfying the semigroup conditions, the equations above define an isolated complete intersection singularity whose link is the $\Z$HS corresponding to that splice diagram.
\end{theorem}

One calls these \emph{singularities of splice type}.  The proof is by induction on the number of nodes; the key step starts with an ``end node'' (all but one of whose emanating edges is a leaf), takes a weighted blow-up, and uses inductively the result for a subdiagram.  Again, setting a variable equal to $0$ gives a monomial curve corresponding to a knot in the link.

How general are singularities of splice type among singularities with $\Z$HS link? By an important general result of Popescu-Pampu \cite{pp}, every splice diagram as in Theorem 2.2 is the link of some \emph{Gorenstein} singularity.    One should keep in mind a class of examples in the paper \cite{nem} of N\'emethi, Luengo, and Melle-Hernandez.
\begin{examples*}
    {\bf a.}~~There exists a Gorenstein singularity, not a complete intersection, whose
    link is the Brieskorn sphere $\Sigma(2,13,31).$

{\bf b.}~~There exists a Gorenstein singularity, not a complete intersection, whose
    link is a homology sphere but which does not satisfy the semigroup
    conditions.
    \end{examples*}
    
An analytic condition distinguishes singularities of splice type from other singularities with the same link.   Each leaf of the splice (or resolution)
diagram gives a knot in $\Sigma$, unique up to isotopy.  A key point
in proving that splice diagram equations
give integral homology sphere links is to show that the variable
$z_{i}$ associated to a leaf cuts out the corresponding knot in
$\Sigma$.  In other words, the curve $C_{i}$ given by $z_{i}=0$ is
irreducible, and its proper transform $D_{i}$ on the minimal good
resolution is smooth and intersects transversely the exceptional curve
corresponding to the leaf of the splice diagram.  Then the
existence of such ``end-curve functions'' implies the semigroup condition on the
splice diagram.

\begin{theorem} \label{th:ends}(\cite{nw2},(8.1))
  Let $(X,0)$ be
  a normal surface singularity whose link $\Sigma$ is an integral
  homology sphere.  Assume that for each of the $t$ leaves $w_i$ of
  the splice diagram $\Delta$ of $\Sigma$, there is a function
  $z_{i}$ inducing the end knot as above.  Then
  \begin{enumerate}
   \item $\Delta$ satisfies the semigroup condition
   \item $X$ is a complete intersection of embedding dimension $\leq t$
   \item $z_{1},\cdots,z_{t}$ generate the maximal ideal of the local
     ring of $X$ at $0$, and $X$ is a complete intersection of splice
     type with respect to these generators.
  \end{enumerate}
\end{theorem}

 Thus, some natural open questions are as follows:  
 \begin{questions*}
 {\bf a.}~~Does there exist a complete intersection singularity for the two-node splice diagram based on $(2,3,1; 2,3,37)$, where the semigroup condition fails?
 
 {\bf b.}~~ Is every complete intersection with $\Z$HS link of splice type?
 
 {\bf c.}~~Assuming there are complete intersections with $\Z$HS link that are not of splice type, do they satisfy the Casson Invariant Conjecture?
 \end{questions*}

 \section{Casson Invariant and Milnor fiber conjectures}
 What about the CIC for these singularities of splice type?  It was realized early on that it is easier to restate it in terms of the geometric genus $p_g$ of the singularity (as opposed to the signature).  Theorem 2 of \cite{nw2} proves the conjecture for  singularities for which the nodes of the splice diagram lie on one line.   The general case was proved by N\'emethi and Okuma, using a more powerful induction.

\begin{theorem}\cite{n-o}  The Casson Invariant Conjecture is true for singularities of splice type.
\end{theorem}

As indicated previously, the validity of the CIC for singularities of splice type suggests that there is a construction of ``splicing Milnor fibers'' for which the additivity of the signatures is geometrically clear.  

Suppose the equations $f_i(z_1,\dots,z_n)=0$, $i=1,\dots,n-2$ define a singularity $X$
of splice type, corresponding to a splice diagram
$\Delta$. Then the curve $z_j=0$ cuts out in $\Sigma$ the knot $K_j$
corresponding to the $j$-th leaf of $\Delta$.  The Milnor
fiber of the singularity at $0$ of the complete intersection curve
$(f_1,\dots,f_{n-2},z_j)^{-1}(0)$ can be considered as a surface $G_j$ in the Milnor fiber of $X$, with $\partial G_j=G_j\cap \Sigma=K_j$.
We proposed in \cite{nw2}, \S6 a conjectural iterative description of a Milnor fiber in terms
of the Milnor fibers of simpler complete intersection surface
singularities and fibers $G_j$ as above lying in their boundaries.

Let $\Delta$ be a splice diagram satisfying the semigroup conditions, and write
$\Sigma=\Sigma_1~\raise5pt\hbox{$\underline{K_1\quad K_2}$}~\Sigma_2$
to represent two homology spheres determined by cutting $\Delta$ at an edge to
form two diagrams with distinguished leaves. Then
$\Delta_1$ and $\Delta_2$ also satisfy the semigroup condition, so
$\Sigma_1$ and $\Sigma_2$ are both complete intersection singularity
links given by equations of splice type.  They thus have
Milnor fibers, say $F_1$ and $F_2$, with $\partial
F_i=\Sigma_i$.

\def\Fconj{\overline F}
For the knot $K_1\subset \Sigma_1$, consider as above $G_1\subset F_1$ , as well as a tubular neighborhood $G_1\times D^2\subset F_1$. Similarly, consider $D^2\times G_2\subset F_2$.



Denote
$$F_1^o:=F_1-(G_1\times( D^2)^o),\quad F_2^o:=F_2-(
(D^2)^o\times G_2),\quad$$
so $\partial F_1^o$ is the union of $G_1\times
S^1$ and the exterior (complement of an open tubular neighborhood) of
the knot $K_1\subset \Sigma_1$, and similarly for $\partial F_2^o$.


\begin{Mconjecture*} The Milnor fiber $F$ is homeomorphic to the result $\Fconj$ of pasting:
  $$\Fconj:=F_1^o\cup_{G_1\times S^1} (G_1\times G_2)\cup_{S^1\times
    G_2}F_2^o,$$
  where we identify $G_1\times S^1$ with
  $G_1\times\partial G_2$ and $S^1\times G_2$ with $\partial G_1\times
  G_2$.

\end{Mconjecture*}

The construction of $\Fconj$ extends splicing on the boundary links.  It is not hard to see that signatures add, specifically $$\text{sign}\  \Fconj=\text{sign}\  F_1+\text{sign}\ F_2.$$   Thus the Milnor Fiber Conjecture implies the Casson Invariant Conjecture.

in \cite{nw2}, Theorem 8.2, the Milnor Fiber Conjecture is verified for the previously discussed singularities $\{z^n+f(x,y)=0\}$ with $\Z$HS link.  In \cite{lam}, the conjecture was proved for iterated suspensions.  Recently, M. A. Cueto, P. Popescu-Pampu, and D. Stepanov \cite{cueto}, \cite{mfc} have announced a proof of the general conjecture, which combines tools from tropical and logarithmic geometry (in the sense of Fontaine and Illusie).

\section{Splice-quotient singularities}

In 2000, as we worked through the basics of singularities of splice type, we began to look for generalizations of the unexpected result of Theorem 1.2, that the UAC of a weighted-homogeneous singularity with $\Q$HS link is a complete intersection. 

From the point of view of classification of
singularities, the easiest examples not covered by Theorem 1.2 are the \emph{quotient-cusps}. These are
singularities whose resolution graphs have the form
$$
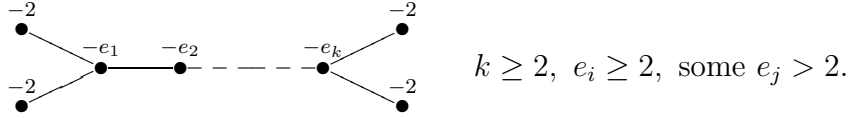
\xymatrix@R=6pt@C=24pt@M=0pt@W=0pt@H=0pt{
  \overtag{\bullet}{-2}{8pt}\lineto[dr] && &&&
  \overtag{\bullet}{-2}{8pt}\lineto[dl]\\
  &\overtag{\bullet}{-e_1}{8pt}\lineto[r]
  &\overtag{\bullet}{-e_2}{8pt}\dashto[r]&\dashto[r]&
  \overtag{\bullet}{-e_k}{8pt}&&k\ge2,~~e_i\ge2, \text{ some }e_j>2.\\
  \overtag{\bullet}{-2}{8pt}\lineto[ur] && &&&
  \overtag{\bullet}{-2}{8pt}\lineto[ul]}$$
They are rational singularities, log-canonical, and \emph{taut}, i.e, the topology of the link
determines the analytic type. 
The above quotient-cusp is double-covered by the cusp singularity
whose resolution graph is
$$
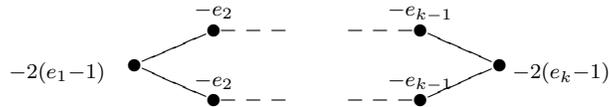
\xymatrix@R=8pt@C=24pt@M=0pt@W=0pt@H=0pt{
  &\overtag{\bullet}{-e_2}{8pt}\dashto[r]
  &&&\overtag{\bullet}{-e_{k-1}}{8pt}\dashto[l]\\
  \lefttag{\bullet}{-2(e_1-1)}{8pt}\lineto[ur]\lineto[dr] &&&&&
  \righttag{\bullet}{-2(e_k-1)}{8pt}\lineto[ul]\lineto[dl]\\
  &\overtag{\bullet}{-e_2}{8pt}\dashto[r]
  &&&\overtag{\bullet}{-e_{k-1}}{8pt}\dashto[l]}$$
It follows that the
universal abelian cover is also a cusp.  

It is easy to
determine, given the resolution graph, when a cusp is a complete
intersection; and, this was the case for the UAC for the few examples we checked.   However, it required a detailed and lengthy analysis of the unimodular matrix classifying the quotient-cusp to prove the following 
\begin{theorem} \cite{nw4} The
  universal abelian cover of a quotient-cusp is a complete
  intersection of embedding dimension $4$.
  \end{theorem} 
Though we wrote down equations and group action for the UAC, at the time (2001) we'd not yet formulated a general method for doing so, and our result seemed ad hoc.  

Given the resolution diagram $\Gamma$ of a singularity with $\Q$HS link $\Sigma$, the method described after Theorem 2.2 allows one to pass from $\Gamma$ to a splice diagram $\Delta$.   Only in the $\Z$HS case can one insure that the weights around a node are relatively prime, and only in that case can one recover $\Gamma$ from $\Delta$.  

Nonetheless, the same approach as above in the $\Z$HS case means that assuming the semigroup conditions on the associated splice diagram, one can write down a set of \emph{splice diagram equations}(\cite{nw1},(2.6)), which give isolated complete intersection surface singularities in $\C^t$, where $t$ is the number of ends of $\Gamma$ or $\Delta$.

When $\Gamma$ is star-shaped, Neumann's Theorem 1.2 proves that these equations give the UAC of the weighted homogeneous singularity of graph $\Gamma$; and, his result describes the group action.  Neglecting this group action in the general $\Q$HS case, the following suggestive assertion appears not to have been proved in full generality:

\begin{conjecture}  If the graph $\Gamma$ satisfies the semigroup conditions, then the \emph{link} of the complete intersection singularity defined by splice diagram equations is the UAC of $\Sigma.$
\end{conjecture}

A difficulty is that splice diagram equations depend only on the splice diagram $\Delta$, while the covering group of the UAC is the discriminant group $D(\Gamma)$, not computable from $\Delta$.  Recall its construction: let $\mathbb E=\oplus_{j=1}^n \Z E_j$ be the lattice spanned by all the exceptional curves, and $\mathbb E^*$ the dual lattice, with dual basis $\{e_i\}$ defined by $e_i(E_j)=\delta_{ij}$.  Then $D(\Gamma)=\mathbb E^*/\mathbb E\cong H_1(\Sigma)$; viewing $\Sigma$ as the boundary of a tubular neighborhood of $E$,  $e_i$ represents an oriented circle over a point of $E_i$.  The pairing of $\mathbb E^*/\mathbb E$ into $\Q/\Z$ is the topological linking pairing on $H_1$.  The key is to construct a natural representation of $D(\Gamma)$.

\begin{proposition} (\cite{nw1},(5.2)) Let $e_{1},\dots,e_{t}$ be the elements of 
  \/ $\E^{\star}$ corresponding to the $t$ leaves of
  $\Gamma$.  
\begin{enumerate}
\item The homomorphism $\E^{\star}\rightarrow \Q^{t}$
  defined by
  $$e\mapsto (e\cdot e_{1},\dots,e\cdot e_{t})
  $$
  induces an injection
  $$D(\Gamma)=\E ^{\star}/\E ~\hookrightarrow~(\Q /\Z )^{t}$$
\item Exponentiating each $\Q/\Z\hookrightarrow \C^*$ via $r \mapsto exp(2\pi i r)$ provides a faithful diagonal representation $D(\Gamma)\hookrightarrow (\C^*)^t$.
\end{enumerate}
\end{proposition}

Now, assume $\Gamma$ is a graph with $t$ ends.  We have a diagonal representation of the discriminant group in $\C^t$. If  $\Delta$ satisfies the semigroup conditions, one has splice type equations in $\C^t$.  What we need is to be able to choose those equations which behave equivariantly with respect to the group action.

\begin{definition}(Congruence conditions) Let $\Gamma$ be a resolution diagram, yielding a splice diagram
  $\Delta$ satisfying the semigroup conditions.  Then $\Gamma$
  satisfies the \emph{congruence conditions} if for each node $v$, one
  can choose for every adjacent edge $e$ an admissible monomial
  $M_{ve}$ so that $D(\Gamma)$ transforms each of these monomials
  according to the same character.
\end{definition}
Clearly, the semigroup plus congruence conditions together mean that the discriminant group acts on appropriate splice diagram complete intersection singularities (higher order terms must also be $D(\Gamma)$-equivariant).  Here is the main result:
\begin{theorem}[Splice-quotient singularities](\cite{nw1},(7.2)) Suppose $\Gamma$ satisfies the semigroup and the
  congruence conditions. Then:
\begin{enumerate}
\item Splice diagram equations for $\Gamma$ define an isolated
complete intersection singularity $(X,0)$.
\item The discriminant group $D(\Gamma)$ acts freely on a
punctured neighborhood of $0$ in $X$.
\item $Y=X/D(\Gamma)$ has an isolated normal surface singularity,
and a good resolution whose associated dual graph is $\Gamma$.
\item $X\rightarrow Y$ is the universal abelian covering.
  \item $X\rightarrow Y$ maps the curve $z_{w}=0$ to an irreducible
  curve, whose proper transform on the good resolution of $Y$ is
  smooth and intersects the exceptional curve transversally, along
  $E_{w}$. In fact the function $z_w^{\disc(\Gamma)}\equiv x_w$, which is\comment{added the
  order of vanishing}
  $D(\Gamma)$--invariant and hence defined on $Y$, vanishes to order
  $\disc(\Gamma)$ on this curve.\label{item:main-endcurve}
\end{enumerate}
\end{theorem}

A singularity $(Y,0)$ which arises in this way is called a \emph{splice-quotient singularity}.

Of course, there are many things to be checked!  The proof is by induction on the number of nodes, and one needs to consider certain non-minimal resolutions.  

We conjectured that rational surface singularities were splice quotients; but the first proof of this fact was due to T. Okuma, who combined the semigroup and congruence conditions in a novel way.

\begin{theorem}\cite{okuma2} Rational surface singularities, and minimally elliptic singularities with $\Q$HS links, are splice quotient singularities.
\end{theorem}

Eventually, we were able to prove a general result, which included Okuma's, identifying the crucial analytic property of a splice-quotient $(Y,0)$: for every end of the graph, there is an ``end-curve function''on $Y$, as in $(5)$ of Theorem 5.4.  This generalized Theorem 3.4 above, but the proof was far more difficult because of the group action.   This is the 

\begin{ECtheorem}\cite{nw3}
  Let $(Y,0)$ be a normal surface singularity with $\Q$HS link
  $\Sigma$.  Suppose that for each leaf $w$ of the resolution diagram
  $\Gamma$ there exists a corresponding end curve function $x_w:Y\rightarrow \C$ which cuts out an end knot $K_w\subset \Sigma$ (or
  end curve) for that leaf. Then $(Y,0)$ is a splice-quotient
  singularity, and a choice of a suitable root $z_w$ of $x_w$ for each
  $w$ gives coordinates for the splice-quotient description.
\end{ECtheorem}

We close with a caveat:  As indicated in Theorem 10.1 of \cite{nw1}, the general splice type equation can be found using any choice of admissible monomials, as long as one allows  the addition of higher order equivariant terms.   This implies that in the $\Z$HS situation, what look like ``equisingular'' deformations of singularities of splice type are also of splice type.   However, even for a weighted homogeneous $\Q$HS singularity, a positive weight deformation of the defining equations might no longer be a splice-quotient.   For example (\cite{nw3},(10.4)), $z^2=x^4+y^9+txy^7$ is a positive weight deformation, but not a deformation of splice-quotients; functions on the resolution lift under deformation, but for $t\neq 0$ they are no longer end-curve functions.

\end{document}